\documentclass[11pt]{gtart}
\usepackage{amsmath, amssymb, graphicx, epsfig, verbatim}
\usepackage{epic}

\newtheorem{thm}{Theorem}[section]

\newtheorem{cor}[thm]{Corollary}
\newtheorem{lem}[thm]{Lemma}
\newtheorem{prop}[thm]{Proposition}

\theoremstyle{definition}

\numberwithin{equation}{section}

\newcommand{\bfn}{{\mathbb {N}}}
\newcommand{\bfz}{{\mathbb {Z}}}

\newcommand{\cpkk}{{\overline {{\mathbb C}{\mathbb P}^2}}}
\newcommand{\cpk}{{\mathbb {CP}}^2}
\newcommand{\cphat}{{\mathbb {CP}}^2\# 6{\overline {{\mathbb C}{\mathbb P}^2}}}

\begin{document}

\title{ Small exotic 4--manifolds with $b_2^+=3$}

\author{Andr\'{a}s I. Stipsicz}
\address{R\'enyi Institute of Mathematics\\
Hungarian Academy of Sciences\\
H-1053 Budapest\\ 
Re\'altanoda utca 13--15, Hungary and\\
Institute for Advanced Study, Princeton, NJ}

\email{stipsicz@renyi.hu, stipsicz@math.ias.edu}

\secondauthor{Zolt\'an Szab\'o}
\secondaddress{Department of Mathematics\\
Princeton University,\\
 Princeton, NJ, 08540}
\secondemail{szabo@math.princeton.edu}

\begin{abstract}
We construct an infinite  family of simply connected,
pairwise nondiffeomorphic  4--manifolds,
all homeomorphic to $3\cpk \# 9 \cpkk$. 
Similar ideas provide examples of 4--manifolds  with $b_2^+=3$,
$b_2^-=8$, vanishing first homology and nontrivial Seiberg--Witten
invariants.
\end{abstract}

\maketitle

\section{Introduction}
Many simply connected, smoothable topological 4--manifolds (with
$b_2^+$ odd) are known to admit infintely many distinct smooth
structures. The existence of such \emph{exotic} structures are,
however, less clear when the Euler characteristic of the 4--manifold
is small, for example if the manifold is homeomorphic to the (blow--up
of the) complex projective plane. Using the rational blow--down
construction \cite{FS1} together with knot surgery \cite{FSknot} in
double node neighborhoods \cite{FSuj}, many new smooth 4--manifolds
have been discovered with $b_2^+=1$ and $b_2^-\geq 5$ \cite{FSuj, P,
PSS, SS}. Similar ideas can be applied to get examples of irreducible
exotic simply connected 4--manifolds with $b_2^+=3$ and relatively
small $b_2^-$.  The study of exotic structures on simply connected
manifolds with $b_2^+=3$ has a rich history. Recall that the $K3$
surface $E(2)$ is such an example with $b_2^-=19$. Applying
appropriate logarithmic transformations on $E(2)$ it was shown that
$3\cpk \# 19 \cpkk$ admits infinitely many smooth structures
\cite{FMbook, MO, SSz}. Later works, relying on Gompf's symplectic
normal connected sum operation, together with Donaldson theory showed
that the topological manifolds $3\cpk \# n \cpkk$ with $n\geq 14$
admit infinitely many smooth structures \cite{Screll, Sztop}. More
recent results of Park \cite{Doug1, Doug2, Doug3} proved the same
statement with $n\geq 10$. Our main result in this paper improves this
bound:

\begin{thm} \label{t:m1}
The simply connected topological 4--manifold $3\cpk \# 9 \cpkk$ admits 
infinitely many distinct smooth structures.
\end{thm}

By modifying the construction used in the proof of Theorem~\ref{t:m1},
we can define a set of 4--manifolds which still have $b_2^+=3$, but
their Euler characteristic is smaller than the above examples. In
these cases, however, we were unable to show that the manifolds are
simply connected.

\begin{thm} \label{t:m2}
There are infinitely many pairwise nondiffeomorphic smooth, closed
4--manifolds with vanishing first homology, $b_2^+=3$, $b_2^-=8$
and nontrivial Seiberg--Witten invariants.
\end{thm}

The proof of the results will involve two steps. First we desribe how
to construct the manifolds claimed and then we use Seiberg--Witten
theory in proving that they are nondiffeomorphic. In the construction
we will apply mapping class group arguments and the theory of
Lefschetz fibrations.  Using the knot surgery construction then we can
identify configurations of curves in the resulting 4--manifolds which
can be rationally blown down, leading us to the desired examples.

\section{Elliptic fibrations and mapping class groups}

Recall that a genus--1 Lefschetz fibration $f\colon X^4\to S^2$ can be
described by the word in the mapping class group $\Gamma _1$ of the
2--torus $T^2$ corresponding to the monodromy presentation of
$f$. More precisely, a point $x\in X$ where $df$ is not onto (called a
singular point of the fibration) gives rise to a singular fiber
$f^{-1}(f(x))$, and the monodromy of the fibration around such a fiber
can be given by the composition of Dehn twists along the circles
corresponding to the vanishing cycles of the singular points of the
fiber. By traversing through the singular fibers in a counterclockwise
manner relative to a fixed base point $s_0\in S^2$, we get the above
mentioned word describing the fibration.  Notice that we do not assume
that $f$ is injective on the set of its singular points, that is, a
singular fiber can contain more than one singular points.  The
assumption that the map $f$ is a Lefschetz fibration implies that the
vanishing cycles corresponding to the singular points in one fixed
singular fiber can be chosen to be disjoint.

It is known that the mapping class group $\Gamma _1$ can be generated by 
two elements $a,b\in \Gamma _1$ which are subject to the two relations
\[ 
aba=bab \qquad {\mbox{ and }} \qquad (ab)^6=1.
\]
In fact, $\Gamma _1$ can be shown to be isomorphic to $SL(2; \bfz )$ 
by mapping $a$ to 
$\left(\begin{smallmatrix}
 1 & 1 \\
 0 & 1 
\end{smallmatrix}\right)$ and $b$ to
$\left(\begin{smallmatrix}
 1 & 0 \\
 -1 & 1 
\end{smallmatrix}\right)$. 
Since the forgetful map from the mapping class group $\Gamma _1 ^1$ of the
2--torus with one marked point to $\Gamma _1$ is an isomorphism, we get that
any genus--1 Lefschetz fibration admits a section. 

Genus--1 Lefschetz fibrations were classified by Moishezon
\cite{Mois}, who showed that after a possible perturbation such a
fibration over $S^2$ is equivalent to one of the fibrations given by
the words $(ab)^{6n}$ ($n\in \bfn$) in $\Gamma _1$. The resulting
4--manifold is usually called $E(n)$ (the simply connected elliptic
surface with section and of holomorphic Euler characteristic $n$), and
$E(2)$ is the famous $K3$ surface.  It can be shown that a section of
$E(n)\to S^2$ has self--intersection $-n$.

Following \cite{HKK} we call a fiber with monodromy conjugate to $a^k$
of \emph{type $I_{k}$} ($k\in \bfn$).  When $k=1$, the corresponding
fiber is also called a \emph{fishtail} fiber.  It is easy to see that
topologically a singular fiber of type $I_k$ ($k\geq 2$) is a plumbing
of $k$ smooth 2--spheres of self--intersection $-2$ plumbed along a
circle (see \cite[page 35]{HKK}), while a fishtail fiber is an
immersed 2--sphere with one positive double point.  Since in $T^2$
nonisotopic simple closed curves necessarily intersect each other, it is
easy to see that a genus--1 Lefschetz fibration can have only
$I_k$--fibers as singular fibers. The fibration $f$ can be perturbed
near a singular fiber $F$ of type $I_k$ into a fibration $f'$ which is
the same as $f$ outside of $F$ but breaks $F$ into two singular fibers
of types $I_{k_1}$ and $I_{k_2}$ with $k_1+k_2=k$.  This fact implies
that a generic genus--1 Lefschetz fibration admits only fishtail
fibers. In our study however, we will find it most helpful to
understand what kind of other singular fibers an elliptic fibration
can admit.

\section{The construction for $b_2^-=9$}

We start with a proposition showing the existence of a particular
genus--1 fibration on $E(2)$.

\begin{prop}\label{p:mcg}
There exists an elliptic Lefschetz fibration on the $K3$ surface
$E(2)$ with a section, a singular fiber $F$ of type $I_{16}$, three
singular fibers $F_1, F_2,F_3$ of type $I_2$ and two further fishtail
fibers.
\end{prop}
\begin{proof}
It is not hard to see that the word $(ab)^{12}$ (defining a genus--1
Lefschetz fibration on the $K3$ surface $E(2)$) in the mapping class
group $\Gamma _1$ of the torus is equivalent to
\[
a^4ba^2b^2a^2b^2a^4ba^2b^2a^2.
\]
(Alternatively, by substituting $a$ and $b$ with the matrices they
correspond to under the map $\Gamma _1 \to SL(2; \bfz )$, we can check
that the above product is equal to the identity matrix.)  By
collecting the powers of $a$ in the front using conjugation, we get
$a^{16}$ followed by the product of three squares of some conjugates
of $b$ and two further conjugates of $b$. Since a conjugate of $b$ by
a word $x\in \Gamma _1$ corresponds to the Dehn twist along the image
under the diffeomorphism $x$ of the curve inducing $b$, the
proposition follows.  As we already mentioned, genus--1 Lefschetz
fibrations always admit sections.
\end{proof}

Perturb first the above fibration near the $I_2$ fibers $F_i$ 
($i=1,2,3$) in a way that these give rise to fishtail fibers
$F_i', F_i ''$ with isotopic vanishing cycles. Let us denote the
resulting Lefschetz fibration by $f\colon E(2) \to S^2$.


Let $K_1,K_2,K_3$ be three twist knots as depicted in \cite{FSuj}.
Let $Z_{K_1,K_2,K_3}$ denote the 4--manifold we get after performing
three knot surgeries with knots $K_1,K_2,K_3$ along three regular
fibers in the fibration on the $K3$ surface found above.  Because of
the existence of fishtail fibers in the complement with nonisotopic
vanishing cycles, we conclude that $\pi _1(Z_{K_1,K_2,K_2})=1$. If we
perform the surgeries in the double node neighborhoods near the
fishtail fibers $F_1',F_1''$, and $F_2',F_2''$ and $F_3',F_3''$
respectively, then we can find a 'pseudo--section' $S$ as in
\cite{FSuj} which is an immersed sphere of self--intersection $(-2)$
with three positive double points, intersecting the further fishtail
fibers and the $I_{16}$ fiber $F$ transversally. Next smooth the
intersections of this pesudo--section $S$ with two further fishtail
fibers. The result is an immersed sphere of self--intersection 2
having 5 positive double points. Now blow up the 4--manifold
$Z_{K_1,K_2,K_3}$ in the five double points of this sphere, and find
an embedded sphere of self--intersection $-18$ in $Z_{K_1,K_2,K_3}\# 5
\cpkk$. Let $C$ denote the tubular neighborhood of the linear plumbing
of spheres given by this $(-18)$--sphere together with 14 of the
$(-2)$--spheres in the $I_{16}$ fiber $F$. It is not hard to see that
$C$ is diffeomorphic to $C_{16,1}$ in the notation of \cite{SS},
cf. also \cite{FS1, Pratb}.  It is then easy to show that $\partial
C=\partial C_{16,1}$ can be given as the oriented boundary of the
rational ball $B_{16,1}$, see \cite{CH, FS1, SS}.  Define
$X_{K_1,K_2,K_3}$ as the rational blow--down of $Z_{K_1,K_2,K_3}\# 5
\cpkk$ along $C$, that is,
\[
X_{K_1,K_2,K_3}=(Z_{K_1,K_2,K_3}\# 5\cpkk -{\mbox { int }}C)\cup B_{16,1}.
\]

\begin{prop}\label{p:hom}
For any twist knots $K_1,K_2,K_3$ the 4--manifold $X_{K_1,K_2,K_3}$ is
homeomorphic to $3\cpk \# 9 \cpkk$.
\end{prop}
\begin{proof}
Notice that $Z_{K_1,K_2,K_3}$ is simply connected, and the complement
of $C$ in $Z_{K_1,K_2,K_3} \# 5 \cpkk$ is simply connected since the
$(-2)$--sphere in $F$ intersecting the last $(-2)$--sphere of the
linear chain $C$ provides a hemisphere which contracts the generator
of the fundamental group of $\partial C$.  Since $\pi _1 (\partial
B_{16,1})$ surjects onto $\pi _1 (B_{16,1})$ under the natural
embedding, Van Kampen's Theorem implies that $X_{K_1,K_2,K_3}$ is
simply connected. Now simple signature and Euler characteristics
computation together with Freedman's Theorem \cite{Fr} verifies the
result.
\end{proof}

For short, let $X_n$ denote the 4--manifold $X_{K_1,K_2,K_3}$ if
$K_1=K_2=K_3=$ the $n$--twist knot $T_n$. The following proposition is
true in a wider generality, we restrict our attention to the special
case $K_1=K_2=K_3$ in order to keep our discussion as simple as possible.
Using the results of \cite{FS1, FSknot} the Seiberg--Witten invariants
of $X_n$ can be easily computed. This computation immediately shows

\begin{thm}\label{t:sw}
There are two cohomology classes $\pm L\in H^2 (X_n ; \bfz )$ with the
property that the Seiberg--Witten function $SW_{X_n}$ of $X_n$ is
equal to zero for all classes $K\neq \pm L$ and $SW_{X_n}(\pm L)=\pm
n^3$.
\end{thm}
\begin{proof}
By \cite[Theorem~1.1]{FSknot} the Seiberg--Witten invariant of
$Z_n=Z_{T_n,T_n,T_n}$ can be computed to be equal to
\[
(n\exp (2T)-(2n-1)+n\exp (-2T))^3
\]
(use the facts that the Seiberg--Witten function of the $K3$ surface
is equal to 1 and the Alexander polynomial of $T_n$ is $\Delta
_{T_n}=nt-(2n-1)+nt^{-1}$). This result, together with the blow--up
formula shows that the five--fold blow--up $Z_n \# 5 \cpkk$ has
exactly two Seiberg--Witten basic classes $\pm K$ which evaluate on
the $(-18)$--sphere of the configuration $C$ as $\pm 16$. Moreover,
the value of the Seiberg--Witten function on these basic classes is
$\pm n^3$.  Now \cite[Theorem~8.5]{FS1} implies that $X_n$ has two
basic classes, on which the value of the Seiberg--Witten function is
equal to $\pm n^3$, verifying the result.
\end{proof}

\begin{cor}\label{c:nondiffeo}
The 4--manifolds $X_n$ are pairwise nondiffeomorphic irreducible smooth
4--manifolds. \qed
\end{cor}

\begin{proof}[Proof of Theorem~\ref{t:m1}]
The 4--manifolds $X_n$ provide an infinite family of smooth
4--manifolds all homeomorphic to $3\cpk \# 9 \cpkk$ according to
Proposition~\ref{p:hom}, and by Corollary~\ref{c:nondiffeo} these
manifolds are pairwise nondiffeomorphic. Therefore the set $\{ X_n
\mid n \in \bfn \}$ provides an infinite collection of distinct
smooth structures on $3\cpk \# 9 \cpkk$, hence the proof is complete.
\end{proof}

\subsection{4--manifolds with $b_2^+=3$ and $b_2^-=8$}
Using a variation of the above procedure, closed 4--manifolds with
$b_2^+=3$ and $b_2^-=8$ can be constructed as follows. Consider the
fibration $f\colon E(2) \to S^2$ found above, containing a singular
fiber of type $I_{16}$ and eight fishtail fibers, out of which three
pairs $F_i', F_i''$ ($i=1,2,3$) have isotopic vanishing cycles.
Proceed as before by doing three knot surgeries along three regular
fibers with twist knots $K_1,K_2,K_3$ and --- using the double node
neighborhoods provided by the fishtail fibers $F_i',F_i''$ ---
identify the 'pseudo--section', which is again an immersed sphere
with homological square $-2$ and has three positive double points. As
before, resolve the two positive intersections of this
pseudo--section with the remaining two fishtail fibers, and find the
immersed sphere with 5 double point and homological square 2. Blow up
the 4--manifold at the double points of the pesudo--section, and
consider the resulting sphere of square $-18$ in
$Z_{K_1,K_2,K_3}$. This sphere intersects the $I_{16}$ fiber $F$
transversally in a unique point $P$ which is on the $(-2)$--sphere
$\Sigma _0 \subset F$. Now $\Sigma _0$ is intersected by two other
spheres in $F$, let $\Sigma _1$ be one of them and denote the
intersection point of $\Sigma _0$ and $\Sigma _1$ by $Q$. Apply 17
infinitely close blow--ups at $Q$. The resulting plumbing manifold
$C_{305,17}$ can be given by the linear plumbing
\[
(-18,-19,-2,-2,-2,-2,-2,-2,-2,-2,-2,-2,-2,-2,-2,-2,-3,-2,-2,
\]
\[
-2,-2,-2,-2,-2,-2, -2,-2,-2,-2,-2,-2,-2,-2).
\] 
As before, it is routine to see that the boundary
$\partial C_{305,17}$ can be given as the boundary of the rational
ball $B_{305,17}$, hence we can blow it down, resulting the 
4--manifold $Y_{K_1,K_2,K_3}$. Since the normal circle of the 
pseudo--section (which circle generates the first homology
of $\partial C_{305,17}$) vanishes in the homology of the
complement $Z_{K_1,K_2,K_3}\#22\cpkk - {\mbox { int }}C_{305,17}$
(as it is shown by a regular fiber), the Mayer--Vietoris sequence
implies that $H_1(Y_{K_1,K_2,K_3}; \bfz )$ vanishes. Simple
Euler characteristic and signature computations now imply
\begin{lem}
The 4--manifolds $Y_{K_1,K_2,K_3}$ have $b_2^+=3$ and $b_2^-=8$. \qed
\end{lem}

Finally, by computing Seiberg--Witten invariants of these manifolds in
the fashion it was done in Theorem~\ref{t:sw}, we get that the
4--manifolds $Y_n=Y_{T_n, T_n,T_n}$ are pairwise nondiffeomorphic, all
with nonvanishing Seiberg--Witten invariants, verifying the claim of
Theorem~\ref{t:m2}.

\end{document}